\documentclass[10pt]{article}
\usepackage{amsmath}
\usepackage{amscd}
\usepackage{amsfonts}
\usepackage{amssymb}
\usepackage{amsthm}
\usepackage{mathrsfs}
\usepackage{latexsym}
\usepackage{graphicx}
\usepackage{times}
\usepackage{natbib}
\usepackage{fullpage}
\usepackage{setspace}
\newcommand{\mbs}[1]{\boldsymbol{#1}}
\title{Double-mixing semiparametric logistic regression with unknown sizes}
\author{Wei Zhang}
\date{}
\DeclareMathOperator*{\argmax}{argmax}
\DeclareMathOperator*{\argmin}{argmin}
\begin{document}
\maketitle
\centerline
{Department of Statistics, University of California, Riverside, CA
, 92521}
\centerline{wxz118@yahoo.com}

\begin{abstract}
Binomial data with unknown sizes 
often appear in biological and medical sciences and  are 
usually overdispersed. All previous methods used parametric 
models and only considered  
overdispersion due to the variation of sizes. 
The proposed semiparametric model
 considers overdispersion due to the variation of sizes and  
that of probabilities. By doing this,  
it can  include variations caused by observations, 
missing covariates, and random measurement errors in covariates.
An Expectation Conditional Maximization  algorithm 
is provided to stabilize the loglikelihood optimization. 
Selecting the number of support points of the mixing 
distributions and the bootstrap methods are also discussed. 
Simulation is done to evaluate the performance of the proposed model. 
Two real examples are used to illustrate the proposed model.
\end{abstract}


\section{Introduction}
The study about  the binomial data with unknown  sizes can be 
dated back to Wadley (1949), which made a fictitious experimental data 
set. In the experiment, fruits were 
infested with fruit-fly larvae and exposed to low temperature with varying 
days. Some of the larvae would die from low temperature. 
The number of fruit-flies that were seen to emerge was 
 counted, but the initial 
larvae number was not known. 

A lot of real binomial data with unknown sizes appeared in the literature. 
For example, 
Morton (1981) presented a data set about the disinfestation of 
wheat. Margolin et al. (1981) studied the effects of 
 quinoline on the
number of revertant colonies of Salmonella strain TA98.  
Elder (1996) investigated the relationship between the survival of 
 V79-473 cells and their times in the high heat. 
Bailer and Piegorsch (2000) 
took the effect of nitrofen on
the offspring of \textit{C. dubia} as an example. 

Let $y_i$ denote a binomial random variable 
 with size $n_i$ unknown and probability $p_i$, and 
$\mbs{x}_i$  denote a vector of covariates of length $\varrho$, 
$i=1$, 2, $\dots$, $r$. 
The issue of interest is to investigate the relationship between 
 the covariates $\mbs{x}_i$ and the probabilities $p_i$. 
If overdispersion exists, there will be three cases: 
overdispersion due to only the variability of $p_i$, due to 
only the variability 
of $n_i$, or due to both of them (Elder et al. 1999). 

The binomial data with unknown sizes are usually approximated by 
Poisson distributions (e.g. Wadley 1949 and Margolin et al. 1981). 
However, such an approximation is not reasonable for moderate sizes or
probabilities (e.g. Elder et al. 1999). Anscombe (1949) considered  
overdispersion due to the $n_i$ and provided 
a parametric model based on the negative binomial distribution. 
 Baker et al. (1980) 
treated $y_i$ in the control group as 
a Poisson random variable  with  mean $m$ and  
that in the treatment 
group with mean $m p_i$, where 
a probit dose-response relationship is assumed. 
 Trajstman (1989) 
modified the method of Baker et al. (1980) to allow a logistic dose-response 
relation and incorporated overdispersion by assuming a scaled Poisson 
variance-mean relationship. Based on Baker et al. (1980), 
Morgan and Smith (1992) used a full negative-
binomial distribution and incorporated extra Poisson variation. 
 Kim and Taylor (1994) and Elder et al. (1999) 
developed a quasi-likelihood approach by regarding
$y_i$ as a binomial random variable.
Kim and Taylor (1994) assumed that $E(n_i)=m_i$
and $\text{var}\,(n_i)=m_i\nu$ with $m_i$ known and $\nu\geqslant 1$
unknown. Elder et al. (1999) estimated $m=E(n_i)$
with $\text{var}\,(n_i)=m(1+\nu m)$ and $\nu \geqslant 0$. 
All previous methods used parametric models and only considered 
overdispersion due to the variation of $n_i$.

We  propose a semiparametric model to incorporate overdispersion 
due to the variability of $n_i$ and $p_i$. The $n_i$ are assumed 
to be Poisson distributed with means from an unspecified mixing distribution.
 With the mean of $n_i$ being a random variable, overdispersion 
due to the variation of $n_i$ is taken into account.
It is assumed that $p_i=p(\eta+\mbs{x}'_i\mbs{\beta})$, where the 
$\eta$ are  further assumed to follow another unspecified mixing distribution.
 With $\eta$ being a random variable, we 
include variations caused by observations, missing covariates, 
and random measurement errors in covariates (Follmann and Lambert 1989). 

Section 2 is the method part, and  includes the proposed model, the 
Expectation Conditional Maximization (ECM)  
algorithm, how to select the number of support points, and the bootstrap 
methods. A simulation study is presented in Section 3. Two real examples
 of bioassay are investigated in Section 4. 

\section{Methods}

\subsection{A model}
In the proposed model, a logistic link is assumed, and its inverse is  
\[p_i=p(\eta+\mbs{x}'_i\mbs{\beta})
=\frac{\exp(\eta+\mbs{x}_i'\mbs{\beta})}
{1+\exp(\eta+\mbs{x}_i'\mbs{\beta})}.
\]
The unknown size $n_i$ is assumed to be a Poisson random variable with 
mean $\xi_i$. It can be easily shown that  
 $y_i|\xi_i,\eta\sim \text{Pois}(\xi_ip(\eta+\mbs{x}'_i\mbs{\beta}))$. 
The nuisance parameters $\xi_i$ and  $\eta$ are further assumed to follow 
mixing distributions $G$ and $H$, respectively. 
Since the parameter of interest $\mbs{\beta}$ 
is in an Euclidean space of dimension $\varrho$, 
a semiparametric model arises 
when $G$ and $H$ are nonparametric.
The density of a single observation $(y,\mbs{x})$ is
\[
f(y;\mbs{x},\mbs{\beta},G,H)=\iint f(y;\mbs{x},\mbs{\beta},
\lambda,\alpha)dG(\lambda)dH(\alpha),
\]
where $f(y;\mbs{x},\mbs{\beta},\lambda,\alpha)$ is a Poisson density 
with mean $\lambda p(\alpha+\mbs{x}'\mbs{\beta})$, i.e.,
\[
f(y;\mbs{x},\mbs{\beta},\lambda, \alpha)=
\exp\{-\lambda p(\alpha+\mbs{x}'\mbs{\beta})\}
\{\lambda p(\alpha+\mbs{x}'\mbs{\beta})\}^y/y!, \quad y=0,1,\dots.
\]
The log likelihood can be written as
\begin{equation}
\label{eq:loglike}
\ell(\mbs{\beta},G,H)=\sum_{i=1}^{r}\log f(y_i;\mbs{x}_i,\mbs{\beta},G,H).
\end{equation}

\subsection{An ECM algorithm}

Since any distribution can be approximated by a discrete distribution, 
we assume  $G$ and $H$ are discrete distributions. First, we 
consider that  $G$ and $H$ have a fixed number of support points 
$K_1$ and $K_2$, which will be allowed to change in Section 2.3.
Let $G=\sum_{j=1}^{K_1}\rho_j\delta(\lambda_j)$ and 
$H=\sum_{m=1}^{K_2}\pi_m\delta(\alpha_m)$,  
where 
$\sum_{j=1}^{K_1}\rho_j=1$, $\sum_{m=1}^{K_2}\pi_m=1$,
$\rho_j\geq 0$, $\pi_m\geq 0$, 
$\delta$ is the indicator function, $\lambda_j\in (0,\infty)$,
 and $\alpha_m\in \mathcal{R}$.
Let $\mbs{\rho}=(\rho_1,\,$$\rho_2,\,\dots,\,\rho_{K_1})'$,
$\mbs{\lambda}=(\lambda_1,\,$$\lambda_2,\,\dots,\,\lambda_{K_1})'$,
 $\mbs{\pi}=(\pi_1,\,$$\pi_2,\,\dots,\,\pi_{K_2})'$, 
$\mbs{\alpha}=(\alpha_1,\,$$\alpha_2,\,\dots,\,\alpha_{K_2})'$,
 and $\mbs{\theta}=$$(\mbs{\beta},$ 
$\mbs{\rho}$, $\mbs{\lambda}$, $\mbs{\pi}$, $\mbs{\alpha})$. The log likelihood becomes 
\begin{equation}
\label{eq:disloglike}
\ell(\mbs{\theta})=\sum_{i=1}^r\log \left\{
\sum_{j=1}^{K_1}\sum_{m=1}^{K_2}\rho_j\pi_mf(y_i;\mbs{x}_i,\mbs{\beta},
\lambda_j,\alpha_m)\right\}.
\end{equation}
Since direct maximization of $\ell(\mbs{\theta})$ in (\ref{eq:disloglike}) 
is extremely difficult,  an EM algorithm may be considered. 
 However, the M-step in the EM algorithm may be computationally  
unreliable because of so many parameters.

The ECM algorithm (Meng and Rubin 1993, 
McLachlan and Peel 2000, p148) is promising.
The ECM algorithm simplifies the M-step by 
replacing the complicated M-step
 with five computationally simpler and stabler CM-steps.
But it also keeps the advantage of the EM algorithm, and  
increases the log likelihood in each iteration
 (Meng and Rubin 1993).

Suppose the missing data are the indicator vectors for the pair $(\mbs{x},y)$, 
$\mbs{z}_1=(z_{11},z_{12},...,z_{1K_1})'$ 
and $\mbs{z}_2=(z_{21},z_{22},...,z_{2K_2})'$, 
 where 
$z_{1j}=1$ means that $\xi=\lambda_j$ and $z_{2m}=1$ means that $\eta=\alpha_m$.
Note that $\mbs{z_1}$ and $\mbs{z_2}$ are independently  
multinomial distributed with size one and 
probability $\mbs{\rho}$ and $\mbs{\pi}$, respectively. 
The complete density for a single datum 
$(\mbs{x},y,\mbs{z}_1,\mbs{z}_2)$ is 
\[
\Pi_{j=1}^{K_1}\Pi_{m=1}^{K_2}\left [\rho_j\pi_mf(y;\mbs{x},\mbs{\beta},\lambda_j,
\alpha_m)
\right ]^{z_{1j}z_{2m}}.
\]
The joint complete log likelihood is
\[
\ell_c(\mbs{\theta})=\sum_{i=1}^r\sum_{j=1}^{K_1}\sum_{m=1}^{K_2}z_{i1j}
z_{i2m}\big\{\log\rho_j+\log\pi_m+\log[f(y_i;\mbs{x}_i,\mbs{\beta},\lambda_j, 
 \alpha_m)]\big\}.
\]
The expected conditional complete log likelihood to be maximized
is 
\[
W(\mbs{\theta};\mbs{\theta}^{(0)})=
E_{\mbs{\theta}^{(0)}}\left\{\ell_c(\mbs{\theta})|y_1,y_2,...,y_r\right\}.
\]

In the E-step, the conditional expectation of $z_{i1j}z_{i2m}$ 
is calculated, 
i.e., for  $i=1,2,\dots,r, j=1,2,\dots,K_1, m=1,2,\dots,K_2$,
\[
e^{(0)}_{ijm}=E_{\mbs{\theta}^{(0)}}\big (z_{i1j}z_{i2m}|y_1,y_2,...,y_r \big)
=\frac{\rho^{(0)}_j\pi^{(0)}_mf(y_i;\mbs{x}_i,\mbs{\beta}^{(0)},\lambda^{(0)}_j, 
\alpha^{(0)}_m)}
{\sum_{h_1=1}^{K_1}\sum_{h_2=1}^{K_2}
\rho^{(0)}_{h_1}\pi^{(0)}_{h_2}f(y_i;\mbs{x}_i,\mbs{\beta}^{(0)},\lambda^{(0)}_{h_1},
\alpha^{(0)}_{h_2})}.
 \]

In the CM-step, the  expected conditional complete log likelihood

\begin{align*}
W(\mbs{\theta};\mbs{\theta}^{(0)})
&=\sum_{i=1}^r\sum_{j=1}^{K_1} \sum_{m=1}^{K_2}
e_{ijm}^{(0)}
\log \rho_j+\sum_{i=1}^r\sum_{j=1}^{K_1}\sum_{m=1}^{K_2} 
e_{ijm}^{(0)}\log \pi_m+\\
&\sum_{i=1}^r\sum_{j=1}^{K_1} \sum_{m=1}^{K_2}
e_{ijm}^{(0)}\log f(y_i;\mbs{x}_i,\mbs{\beta},\lambda_j,\alpha_m)\\
&=\text{constant}+\underbrace{\sum_{i=1}^r\sum_{j=1}^{K_1}
\sum_{m=1}^{K_2} e_{ijm}^{(0)}
\log \rho_j}_{T_1(\mbs{\rho})}+\underbrace{\sum_{i=1}^r\sum_{j=1}^{K_1}
\sum_{m=1}^{K_2} e_{ijm}^{(0)}
\log \pi_m}_{T_2(\mbs{\pi})}+\\
&\underbrace{\sum_{i=1}^r\sum_{j=1}^{K_1}\sum_{m=1}^{K_2}
e^{(0)}_{ijm}
\left\{y_i\log\lambda_j+y_i\log p(\alpha_m+\mbs{x}_i'
\mbs{\beta})-\lambda_j p(\alpha_m+\mbs{x}_i'
\mbs{\beta})\right\}}_{T_3(\mbs{\lambda},\,\mbs{\alpha},\,\mbs{\beta})}
\end{align*}
 is maximized 
over $\mbs{\rho}$, $\mbs{\pi}$, $\mbs{\lambda}$, $\mbs{\alpha}$ and $\mbs{\beta}$ sequentially.
Because $\mbs{\rho}$, $\mbs{\pi}$ and $(\mbs{\lambda},\,\mbs{\alpha},\,\mbs{\beta})$
are in $T_1(\mbs{\rho})$, $T_2(\mbs{\pi})$ 
and $T_3(\mbs{\lambda},\,\mbs{\alpha},\,\mbs{\beta})$
 separately, 
 their maximum likelihood estimators (MLEs)  can be found individually. 
The MLE for $\mbs{\rho}$ is
 \begin{equation}
\label{eq:rho}
\rho_j^{(1)}=r^{-1}\sum_{i=1}^r\sum_{m=1}^{K_2}e^{(0)}_{ijm}, \,j=1,2,...,K_1.
\end{equation}
The MLE for $\mbs{\pi}$ is
 \begin{equation}
\label{eq:pi}
\pi_m^{(1)}=r^{-1}\sum_{i=1}^r\sum_{j=1}^{K_1}e^{(0)}_{ijm}, \,m=1,2,...,K_2.
\end{equation}

We will maximize $T_3(\mbs{\lambda},\,\mbs{\alpha},\,\mbs{\beta})$
 over $\mbs{\lambda}$, $\mbs{\alpha}$ and $\mbs{\beta}$ sequentially. 
The conditional MLE for $\mbs{\lambda}$
given $\mbs{\alpha}=\mbs{\alpha}^{(0)}$ and $\mbs{\beta}=\mbs{\beta}^{(0)}$ 
is
\begin{equation}
\label{eq:lambda}
\lambda_j^{(1)}=\frac{\sum_{i=1}^r\sum_{m=1}^{K_2}e^{(0)}_{ijm}y_i}
{\sum_{i=1}^r\sum_{m=1}^{K_2}e^{(0)}_{ijm}
p(\alpha_m^{(0)}+\mbs{x}_i'\mbs{\beta}^{(0)})},\,
j=1,2,...,K_1.
\end{equation}

There are no simple analytic forms for the conditional MLEs 
of $\mbs{\alpha}$ and  $\mbs{\beta}$. 
The conditional MLE for $\mbs{\alpha}$ 
given $\mbs{\lambda}=\mbs{\lambda}^{(1)}$ and 
$\mbs{\beta}=\mbs{\beta}^{(0)}$ is
\begin{equation}
\label{eq:alpha}
\alpha^{(1)}_m=\argmax_{\alpha_m\in \mathcal{R}} 
\sum_{i=1}^r\sum_{j=1}^{K_1}
e^{(0)}_{ijm}
\left\{y_i\log\lambda_j+y_i\log p(\alpha_m+\mbs{x}_i'
\mbs{\beta})-\lambda_j p(\alpha_m+\mbs{x}_i'
\mbs{\beta})\right\}, 
\end{equation}
for $m=1,2,\dots,K_2.$ 
The conditional MLE for $\mbs{\beta}$
given $\mbs{\lambda}=\mbs{\lambda}^{(1)}$ and 
$\mbs{\alpha}=\mbs{\alpha}^{(1)}$
is
\begin{equation}
\label{eq:beta}
\mbs{\beta}^{(1)}=\argmax_{\mbs{\beta}\in \mathcal{R}^\varrho} 
T_3(\mbs{\lambda}^{(1)},\, \mbs{\alpha}^{(1)},\,\mbs{\beta}).
\end{equation}
The function \textit{optim} in R can be used to get 
the MLEs of $\mbs{\alpha}$
 and $\mbs{\beta}$ in equations (\ref{eq:alpha}) and 
(\ref{eq:beta}).

\subsection{Selecting the number of support points}

The maximized
 log likelihood $\ell(\hat{\mbs{\theta}})$ can be increased 
by increasing the number of support points of $G$ or $H$. 
We propose to choose the number of support points by 
minimizing the BIC (e.g., Wang et al. 1996) to obtain 
a reasonable and parsimonious fit to the data, i.e., 
\[
(\widehat{K}_1,\widehat{K}_2)=\argmin_{(K_1,K_2)\in \{1,2,\dots\}^2} 
\{-2\ell (\hat{\mbs{\theta}})+\log(r)\left[2(K_1+K_2)-2+\varrho\right]\}.
\]

Forward model selection is used. For a fixed $K_1$, 
if the BIC stops to decrease for larger $K_2$, the models with 
greater $K_2$ will not be considered for this $K_1$. 
The strategy is the same for $K_2$ fixed and $K_1$ changed. 
From all the models considered, the one with the minimum BIC is chosen.
 
\subsection{The bootstrap method}

The confidence intervals for the regression coefficients $\mbs{\beta}$
can be got by the bootstrap method. 
 The nonparametric bootstrap 
method may be applied for a random design, 
in which one can sample the pairs $(y_i,\mbs{x}_i)$.
For a fixed design, a parametric bootstrap method is recommended.
A resample of size $r$ is generated as follows,
\[
y_i^*\sim f(y;\mbs{x}_i,\hat{\mbs{\beta}},\lambda_i,\alpha_i), 
i=1, 2, \dots, r,
\]
where $ \lambda_i$ and $\alpha_i$ are  random variables drawn from
the estimated mixing distributions $\widehat{G}$ and $\widehat{H}$, 
respectively, where
\[
\widehat{G}=\sum_{j=1}^{\widehat{K}_1}\hat{\rho}_j\delta(\hat{\lambda}_j)
\quad \text{and} \quad
\widehat{H}=\sum_{m=1}^{\widehat{K}_2}\hat{\pi}_m\delta(\hat{\alpha}_m).
\]
\section{Simulation}

In the simulation, 
the parameter setting takes a $2^3$ design, 
\begin{equation*}
\underbrace{\{-2,3\}}_{\beta}
\times \underbrace{\{G_1,G_2\}}_{G}
\times \underbrace{\{H_1,H_2\}}_{H},
\end{equation*}
where $G_1=0.1\delta(100)+0.8\delta(200)+0.1\delta(300)$,  
$G_2=0.5\delta(10)+0.5\delta(50)$, $H_1=0.3\delta(-2)$$+
0.3\delta(0.4)+0.4\delta(3)$, and 
$H_2=0.25\delta(-2)+0.75\delta(1.5)$. 
There is a single covariate $x$ in the simulation. 
For each integer $x$ in $[-5,5]$ , $10$ values of $y$ are drawn independently 
from a Possion distribution with mean $\lambda p(\alpha+x\beta)$,
where $\lambda$ and $\alpha$ are random variables drawn from 
$G$ and $H$. So the sample size $r$ is $110$.

For each parameter setting, $200$ samples are drawn.
Table \ref{tab:simres} presents the simulation results. 
The bias, standard deviation 
and mean square error of $\beta$ are small, and
each $\beta$ falls into its $95\%$ quantile interval, with 
$2.5\%$ and $97.5\%$ quantiles as endpoints.

\begin{table}[ht]
\caption{Simulation results: 
sd stands for standard deviation, qi for $95\%$ quantile 
intervals and mse for mean square error.}
\label{tab:simres}
\begin{center}
\begin{tabular}{crccrccc}
\hline
setting & beta & G & H & bias & sd & qi& mse \\
\hline
1 & $-$2 &$G_1$&$H_1$ &$-$0.00 & 0.08 & ($-$2.15, $-$1.86) & 0.01 \\
2 & $-$2 & $G_1$&$H_2$&$-$0.01 & 0.09 & ($-$2.17, $-$1.85) & 0.01 \\
3 & $-$2 &$G_2$&$H_1$& $-$0.11 & 0.32 & ($-$2.77, $-$1.62) & 0.12 \\
4 & $-$2 &$G_2$&$H_2$& $-$0.07 & 0.24 & ($-$2.69, $-$1.68) & 0.06 \\
5 & 3&$G_1$&$H_1$ & 0.00 & 0.16 & (2.72,  3.27) & 0.02 \\
6 & 3& $G_1$&$H_2$ & 0.02 & 0.16 & (2.72, 3.33) & 0.02 \\
7 & 3 &$G_2$&$H_1$ & 0.17 & 0.68 & (2.20, 4.50) & 0.49 \\
8 & 3 &$G_2$&$H_2$ & 0.03 & 0.45 & (2.36, 4.23) & 0.20 \\
\hline
\end{tabular}
\end{center}
\end{table}

\section{Example}
\subsection{\textit{M. Bovis} data}

Table \ref{tb:bovisdata} is part of Table $1$ of
Trajstman (1989), and also appeared in Morgan and Smith (1992). 
One of the decontaminants, HPC or oxalic acid,  
with a specific concentration was applied on a group of 
 \textit{M. bovis} cells, which were then placed on the culture 
plates for colony formation. After $12$ weeks (at stationarity),
 the number of \textit{M. bovis} colonies were counted,
which is equal to the number of surviving \textit{M. bovis} cells.

\begin{table}
\small
\caption{The \textit{M. bovis} cell survival data.}
\label{tb:bovisdata}
\begin{center}
\begin{tabular}{lccccccccccr}
\hline
 $\%$weight/volume & \multicolumn{10}{c}{No. of \textit{M. bovis} colonies at 
stationarity}
 &sample mean  \\
\hline
 &  &  &  & &  &  &  & &  & &  \\
&\multicolumn{10}{c}{control experiment (no decontaminant)}&\\
  & 52 & 80 & 55 & 50 & 58 & 50 & 43 & 50 & 53 & 54 & 51.8  \\
  & 44 & 51 & 34 & 37 & 46 & 56 & 64 & 51 & 67 & 40 &   \\
  &  &  &  & &  &  &  & &  & &  \\
$\text{[HPC]}$&\multicolumn{10}{c}{decontaminant: HPC}&\\
 0.75 & 2 & 4 & 8 & 9 & 10 & 1 & 0 & 5 & 14 & 7 & 6.0  \\
 0.375 & 11 & 12 & 13 & 12 & 11 & 13 & 17 & 16 & 21 & 2 & 12.8 \\
 0.1875 & 16 & 6 & 20 & 23 & 23 & 39 & 18 & 23 & 33 & 21 & 22.2  \\
 0.09375 & 33 & 46 & 42 & 18 & 35 & 20 & 19 & 29 & 41 & 36 & 31.9  \\
 0.075 & 30 & 30 & 27 & 53 & 51 & 39 & 31 & 36 & 38 & 22 & 35.7  \\
 0.0075 & 53 & 62 & 38 & 54 & 54 & 38 & 46 & 58 & 54 & 57 & 51.4 \\
 0.00075 & 3 & 42 & 45 & 49 & 32 & 39 & 40 & 34 & 45 & 51 & 38.0  \\
 &  &  &  & &  &  &  & &  & &  \\
$\text[Oxalic\,\, acid]$
&\multicolumn{10}{c}{decontaminant: oxalic acid}&\\
 5 & 14 & 15 & 6 & 13 & 4 & 1 & 9 & 6 & 12 & 13 & 9.3  \\
 0.5 & 27 & 33 & 31 & 30 & 26 & 41 & 33 & 40 & 31 & 20 & 31.2 \\
 0.05 & 33 & 26 & 32 & 24 & 30 & 52 & 28 & 28 & 26 & 22 & 30.1 \\
 0.005 & 36 & 54 & 31 & 37 & 50 & 73 & 44 & 50 & 37 &  & 45.8  \\
\hline
\end{tabular}
\end{center}
\end{table}

\begin{table}[ht]
\caption{The BIC for the estimated mixing distributions.}
\label{tb:BIC}
\begin{center}
\begin{tabular}{rrrrr}
\hline
$K_1\backslash K_2$ & 1 & 2 & 3 & 4 \\
\hline
1 & 1061.1 & 991.3 & 977.6 & 985.1 \\
2 & 998.0 & \fbox{976.9} & 978.3 & 988.9 \\
3 & 977.0 & 978.7 & 994.3 &  \\
4 & 984.2 & 988.5 &  &  \\
\hline
\end{tabular}
\end{center}
\end{table}

An ANOVA model is fitted with a separate factor for each level 
of the decontaminants. 
Let $x_j$ denote a factor for the concentration level $j$ of the  
decontaminants. 
It is assumed that the $p_i$ satisfy that
\begin{equation}
\label{eq:pbovis}
\log\Bigl\{\frac{p_i}{1-p_i}\Bigr\}=\eta+\sum_{j=1}^{11}\beta_jx_{ij},
\,\,\, i=1,2,\dots,129,
\end{equation}
where $\eta$ is the control effect and $\beta_j$ is the effect 
 difference between  dose $j$ and the control one, 
$j=1,2,\dots,11$.

From Table \ref{tb:BIC}, the case in which $\widehat{G}$ and $\widehat{H}$ 
with 2 support points has the minimum BIC. Therefore, the MLEs
 of $G$ and $H$ are 
$\widehat{G}=0.04\delta(12.41)+0.96\delta(99.32)$, and 
$\widehat{H}=0.82\delta(-0.03)+0.18\delta(0.63)$, respectively. 

\begin{figure}
\centering
\includegraphics[angle=270,width=12cm]{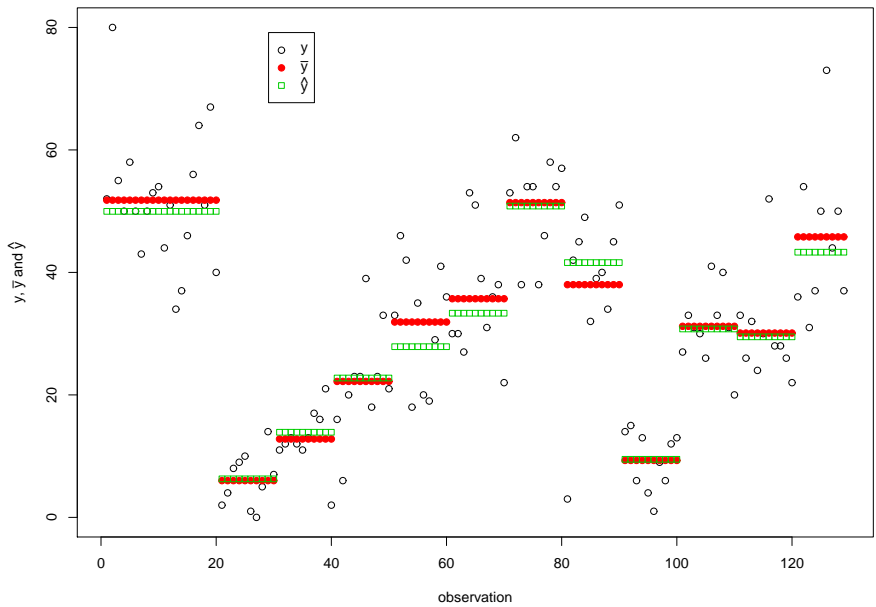}
\caption{The response $y$, sample mean $\bar{y}$ and fitted value $\hat{y}$.}
\label{fig:bovisfit}
\end{figure}

Table \ref{tb:bovisbeta} presents the 
estimated regression coefficients, their bootstrap  
standard errors and $95\%$ confidence intervals from
 $200$ bootstrap resamples.  
The MLEs $\hat{\beta}_6$ and $\hat{\beta}_9$ violate
the monotonic dose-response relationship, i.e.,
 the larger dose do not produce stronger effects here. This
 finding is consistent with the monotonicity violation of their 
sample means in Table \ref{tb:bovisdata}. Nine out of eleven
 confidence intervals 
do not include $0$, so the corresponding doses have  significantly 
stronger negative effects than the control one. 
Two exceptions are those of 
$\beta_6$ and  $\beta_{11}$. An explanation is that 
the  Oxalic acid dose $0.005$  is too small to take any 
different effect on 
the \textit{M. Bovis} cells from the control one. 
The estimates  of $\mbs{\beta}$ in Table \ref{tb:bovisbeta} 
can not be compared to those 
 of Trajstman (1989) and Morgan and Smith (1992), because they used a 
simple linear model in (\ref{eq:pbovis}). 
Figure \ref{fig:bovisfit} presents the responses $y$, the 
sample means $\bar{y}$ and the fitted values $\hat{y}$. 
The model seems fit well.

\begin{table}
\caption{The estimated regression coefficients,
 bootstrapped standard error, and $95\%$ confidence interval 
for the \textit{M. Bovis} data.}
\label{tb:bovisbeta}
\begin{center}
\begin{tabular}{lcrrc}
\hline
dose& $\beta$&MLE & se & 95$\%$ ci\\
\hline
$\textit{HPC}$&&&&\\
 0.75 &$\beta_1$&  $-$2.74 & 0.55 & ($-$4.29, $-$2.36) \\
0.375 &$\beta_2$&  $-$1.86 & 0.51 & ($-$3.38, $-$1.55) \\
0.1875 &$\beta_3$&  $-$1.25 & 0.48 & ($-$2.59, $-$0.96) \\
0.09375 &$\beta_4$&  $-$0.98 & 0.45 & ($-$2.12, $-$0.69) \\
 0.075 & $\beta_5$& $-$0.72 & 0.42 & ($-$1.72, $-$0.48) \\
 0.0075 &$\beta_6$&  0.04 & 0.94 & ($-$0.25,\enspace\, 1.38) \\
 0.00075 &$\beta_7$&  $-$0.35 & 0.38 & ($-$1.14, $-$0.07) \\
$\textit{Oxalic acid}$&&&&\\
5 & $\beta_8$& $-$2.30 & 0.53 & ($-$3.63, $-$1.96) \\
 0.5 &$\beta_9$&  $-$0.84 & 0.45 & ($-$1.99, $-$0.55) \\
 0.05 &$\beta_{10}$&  $-$0.90 & 0.41 & ($-$2.11, $-$0.68) \\
 0.005 &$\beta_{11}$&  $-$0.28 & 0.38 & ($-$0.97, \enspace\, 0.08) \\
\hline
\end{tabular}
\end{center}
\end{table}

\subsection{Jejunal crypt data}

The jejunal crypt data are referred to Table 1 of Elder et al. (1999), which 
are also studied by Kim and Taylor (1994). 
There are $126$ live mice divided into groups, not all of equal sizes. 
The treatment consists of exposing each group of mice to 
 a certain dose of gamma rays, and then 
killing them  to find out the number of surviving crypts. 
The total number of crypts in each mouse is unknown,
because the experiment needs live mice. 
It is assumed that the surviving probabilities $p_i$ satisfy that
\[
\log\left\{\frac{p_i}{1-p_i}\right\}=
\eta+\beta x_i,\quad i=1,2,\dots,126,
\]
 where $x_i$ is the gamma dose.

The BIC are 724.2 for $K_1=K_2=1$, 733.9 for $K_1=2,$ $K2=1$ and 
$K_1=1$, $K_2=2$. Thus, the estimated $\widehat{G}$ is degenerated 
at $\hat{\alpha}=6.705$,
and $\widehat{H}$ at $\hat{\lambda}=196.1$. 
We draw $200$ bootstrap resamples. 
Table \ref{tb:jres} presents the estimation results of 
the proposed model and the previous ones. The bootstrap standard error 
of $\beta$ is quite small and its  $95\%$ confidence interval 
is $(-1.225, -1.029)$. 
All listed estimates of $\beta$ 
 lie in our confidence interval. Since this interval does not 
cover $0$, the $\beta$ in the proposed 
model is significant at the significance level of $0.05$.

\begin{table}
\caption{Jejunal crypt data results from the proposed and
previous approaches 
(logistic regression and Kim's method fix $n_i$ and $E(n_i)$ at 160, 
 respectively; Kim's and Elder's quasi-likelihood method of moments estimates
come from Elder et al. (1999)).}
\label{tb:jres}
\begin{center}
\begin{tabular}{c|ccccc}
\hline
& \multicolumn{4}{c}{estimate (standard error)}\\
\cline{2-5}
&logistic &Kim's&
 Elder's&proposed\\
$\alpha$ & 7.432 (0.175)&7.410 (0.191) &6.727 (0.725)&6.705
\\
\hline
$\beta$ & $-$1.185 (0.024) &$-$1.183 (0.026)& $-$1.126 (0.061)&$-$
1.124 (0.044)\\
$\lambda$ &--- &--- & 194.7 (43.4)&196.1\\
\hline
\end{tabular}
\end{center}
\end{table}

\section{Discussion}

We propose a flexible semiparametric model to incorporate 
 overdispersion due to the variation of $n_i$ and $p_i$. 
The regression coefficients are estimated with the nuisance parameters, the 
mixing distributions in a seamless fashion. 
Although a logistic dose-response relation is assumed, 
it can be extended to other links very easily. When one runs the ECM algorithm, 
a Poisson regression is suggested to run first to obtain 
a good initial value of $\mbs{\beta}$.
 
\nocite{*}
\bibliography{ref}

\begin{thebibliography}{}

\bibitem[Anscombe, 1949]{Anscombe:1949}
Anscombe, F.~J. (1949).
\newblock Note on a problem in probit analysis.
\newblock {\em Annals of Applied Biology}, 36:203--205.

\bibitem[Bailer and Piegorsch, 2000]{Bailer:Piegorsch:2000}
Bailer, A.~J. and Piegorsch, W.~W. (2000).
\newblock From quantal counts to mechanisms and systems: the past, present, and
  future of biometrics in environmental toxicology.
\newblock {\em Biometrics}, 56:327--336.

\bibitem[Baker et~al., 1980]{Baker:1980}
Baker, R.~J., Pierce, C.~B., and Pierce, J.~M. (1980).
\newblock {W}adley's problem with controls.
\newblock {\em GLIM Newsletter}, 3:32--35.

\bibitem[Elder, 1996]{Elder:1996}
Elder, J.~A. (1996).
\newblock {\em Development of quasi-likelihood techniques for the analysis of
  pseudo-proportional data}.
\newblock Unpublished doctoral dissertation, Virginia Commonwealth University,
  Medical College of Virginia, Department of Biostatistics.

\bibitem[Elder et~al., 1999]{Elder:etal:1999}
Elder, J.~A., Carter, W.~H., Gennings, C., and Elswick, R.~K. (1999).
\newblock A quasi-likelihood approach for overdispersed binomial data when
  ${N}$ is unobserved.
\newblock {\em Journal of Agricultural, Biological, and Environmental
  Statistics}, 4:102--115.

\bibitem[Follmann and Lambert, 1989]{Follmann:1989}
Follmann, D.~A. and Lambert, D. (1989).
\newblock Generalizing logistic regression by nonparametric mixing.
\newblock {\em Journal of the American Statistical Association}, 84:295--300.

\bibitem[Kim and Taylor, 1994]{Kim:Taylor:1994}
Kim, D.~K. and Taylor, J. M.~G. (1994).
\newblock Transform-both-sides approach for overdispersed binomial data when
  ${N}$ is unobserved.
\newblock {\em Journal of the American Statistical Association}, 89:833--845.

\bibitem[Margolin et~al., 1981]{Margolin:etal:1981}
Margolin, B.~H., Kaplan, N., and Zeiger, E. (1981).
\newblock Statistical analysis of the {A}mes salmonella/microsome test.
\newblock {\em Proceedings of the National Academy of Sciences}, 78:3779--3783.

\bibitem[McLachlan and Peel, 2000]{McLachlan:2000}
McLachlan, G. and Peel, D. (2000).
\newblock {\em Finite Mixture Models}.
\newblock Wiley.

\bibitem[Meng and Rubin, 1993]{Meng:Rubin:1993}
Meng, X.~L. and Rubin, D.~B. (1993).
\newblock Maximum likelihood estimation via the {ECM} algorithm: a general
  framework.
\newblock {\em Biometrika}, 80:267--278.

\bibitem[Morgan and Smith, 1992]{Morgan:Smith:1992}
Morgan, B. J.~T. and Smith, D.~M. (1992).
\newblock A note on {W}adley's problem with overdispersion.
\newblock {\em Applied Statistics}, 41:349--354.

\bibitem[Morton, 1981]{Morton}
Morton, R. (1981).
\newblock Generalized spearman estimators of relative dose.
\newblock {\em Biometrics}, 37:223--233.

\bibitem[Trajstman, 1989]{Trajstman:1989}
Trajstman, A.~C. (1989).
\newblock Indices for comparing decontaminants when data come from
  dose-response survival and contamination experiments.
\newblock {\em Applied Statistics}, 38:481--494.

\bibitem[Wadley, 1949]{Wadley:1949}
Wadley, F.~M. (1949).
\newblock Dosage-mortality correlation with number treated estimated from a
  parallel sample.
\newblock {\em Annals of Applied Biology}, 36:196--202.

\bibitem[Wang et~al., 1996]{Wang:etal:1996}
Wang, P., Puterman, M.~L., Cockburn, I., and Le, N.~D. (1996).
\newblock Mixed poisson regression models with covariate dependent rates.
\newblock {\em Biometrics}, 52:381--400.

\end{thebibliography}
\bibliographystyle{apalike}

\end{document}